\documentclass[a4paper]{amsart}

\usepackage{amssymb,amsmath,graphics,verbatim}
\usepackage{latexsym}
\usepackage{eucal}

\makeatletter
\@addtoreset{equation}{section}\makeatother

\newtheorem{theo}{Theorem}[section]
\newtheorem{lem}[theo]{Lemma}
\newtheorem{prop}[theo]{Proposition}
\newtheorem{cor}[theo]{Corollary}

\newcommand{\mc}{\mathcal}
\newcommand{\rr}{\mathbb{R}}
\newcommand{\nn}{\mathbb{N}}
\newcommand{\cc}{\mathbb{C}}

\newcommand{\la}{\lambda}
\newcommand{\eps}{\epsilon}

\newcommand{\pl}{\partial}
\newcommand{\p}{\partial}
\newcommand{\x}{\times}

\newcommand{\til}{\widetilde}

\newcommand{\demi}{\frac{1}{2}}
\newcommand{\ndemi}{\frac{n}{2}}
\newcommand{\tra}{\textrm{Tr}}

\newcommand{\ric}{\textrm{Ric}}

\def\qed{\hfill$\square$}
\begin{document}
\title{Inverse problems for Einstein manifolds}
\author{Colin Guillarmou}
\address{Laboratoire J.-A. Dieudonn\'e\\
U.M.R. 6621 du C.N.R.S.\\
Universit\'e de Nice\\
Parc Valrose\\
06108 Nice Cedex 02\\
France}
\email{cguillar@math.unice.fr}

\author{Ant\^onio S\'a Barreto}
\address{Department of mathematics\\
         Purdue University\\
        150 N. University Street, West-Lafayette\\    
         IN 47907, USA}
     \email{sabarre@math.purdue.edu}

%
\begin{abstract}

We show that  the Dirichlet-to-Neumann operator 
of the Laplacian on an open subset of the boundary of a connected 
compact Einstein manifold with boundary determines the manifold up to isometries.  Similarly, for connected conformally compact
Einstein manifolds of even dimension $n+1,$ we prove that  the scattering matrix at energy $n$ 
on an open subset of its boundary  determines the manifold up to isometries.

\end{abstract}
\maketitle
\section{Introduction}\label{intro}
The purpose of this note is to prove two results: first that compact connected Einstein manifolds with boundary are determined modulo isometries from the Dirichlet-to-Neumann map on an open subset of its boundary.    Secondly, that a conformally compact connected Einstein manifolds of even dimension $n+1$ is determined, modulo isometries,  by the scattering matrix  on an open subset  of the boundary. 

The Dirichlet-to-Neumann (DN in short) map $\mc{N}: C^\infty(\pl\bar{X})\to C^{\infty}(\pl\bar{X})$ for the Laplacian on a Riemannian manifold with boundary $(\bar{X},g)$ is defined by solving the Dirichlet problem
\begin{equation}\label{solving}
\Delta_{g}u=0,\quad u|_{\pl\bar{X}}=f 
\end{equation}
where $f\in C^{\infty}(\pl\bar{X})$ is given, then $\mc{N}f:=-\pl_{n}u|_M$ where $\pl_n$ is the interior pointing 
normal vector field to the boundary for the metric $g$. It is an elliptic pseudo-differential operator of order $1$ on the boundary, see for example \cite{LeU}.
Mathematically, it is of interest to know what this map determines about the geometry of the manifold, but $\mc{N}$ can also be interpreted as a boundary measurement of current flux in terms of voltage in electrical impedance tomography.
We refer to \cite{U} for a survey in the field, and to \cite{LUT,LU,LeU,N,SU} for significant results about that problem.\\

Our first result answers a conjecture of Lassas and Uhlmann \cite{LU}
\begin{theo}\label{th1}
Let $(\bar{X}_1,g_1)$ and $(\bar{X}_2,g_2)$ be two smooth connected compact manifolds with respective boundaries 
$\pl \bar{X}_1$ and $\pl \bar{X}_2$. We suppose that $g_1$ and $g_2$ are Einstein with the same constant $\la\in\rr$, i.e. $\ric(g_i)=\la g_i$ for $i=1,2$. Assume that $\pl\bar{X_1}$ and $\pl\bar{X}_2$ contain a common open set $\Gamma$ such that the identity
map ${\rm Id}:\Gamma\subset\bar{X}_1\to \Gamma\subset \bar{X}_2$ is a smooth diffeomorphism. 
If the Dirichlet-to-Neumann map $\mc{N}_i$ of 
$\Delta_{g_i}$ on $\bar{X}_i$ for $i=1,2$ satisfy $(\mc{N}_1f)|_\Gamma=(\mc{N}_2f)|_{\Gamma}$ for any $f\in C_0^\infty(\Gamma)$, then there exists a diffeomorphism $J:\bar{X}_1\to \bar{X}_2$, such that 
 $J^{*}g_2=g_1$. 
\end{theo}

Then we consider a class of non-compact complete Einstein manifolds, but conformal to a compact manifold.   In this case we say that $(X,g)$ is Einstein, with $\operatorname{dim} X=n+1,$ if
$$\ric(g)=-n g.$$ 
 We say 
that a Riemannian manifold $(X,g)$ is conformally compact if $X$ compactifies into a smooth manifold with boundary
$\bar{X}$ and for any smooth boundary defining function $\rho$ of $\bar{X}$, $\bar{g}:=\rho^2g$ extends
to $\bar{X}$ as a smooth metric. Such a metric $g$ is necessarily complete on $X$ 
and its sectional curvatures are pinched negatively outside a compact set of $X.$ 
If in addition the sectional curvatures of $g$ tends to $-1$ at the boundary, we say that $(X,g)$  
is \emph{asymptotically hyperbolic}. 

It has been shown in  \cite{GR,GRL} that if $(X,g)$ is asymptotically hyperbolic, or in particular if
$(X,g)$ is Einstein, then there exists a family of boundary defining functions
$\rho$ (i.e. $\pl\bar{X}=\{\rho=0\}$ and $d\rho|_{\pl\bar{X}}$ does not vanish) such that $|d\rho|_{\rho^2g}=1$
near the boundary. These will be called \emph{geodesic boundary defining functions}. 
Note that, in this case, a DN map can not be defined as in \eqref{solving} since $\Delta_g$ is not an 
elliptic operator at the boundary. 
The natural analogue of the DN map on a conformally compact Einstein manifold
$(X,g)$ is related to scattering theory, at least in the point of view of Melrose \cite{Me}.
We consider an $n+1$-dimensional conformally compact Einstein manifold $(X,g)$ with $n+1$ even. 
Following \cite{GRZ,JSB}, the \emph{scattering matrix or scattering map} in this case, and more generally for asymptotically hyperbolic manifolds, is an operator $\mc{S}: C^{\infty}(\pl\bar{X})\to C^{\infty}(\pl\bar{X})$, 
constructed by solving a Dirichlet problem in a way similar to \eqref{solving}.  This will be discussed in details in section 4.  We show that
for all $f\in C^{\infty}(\pl\bar{X})$, there exists a unique function $u\in C^{\infty}( \bar{X})$ 
such that 
\begin{gather}
\begin{gathered}
\Delta_gu=0  \text{ and } u|_{\pl X}=f.
\end{gathered}\label{poisson}
\end{gather}
Since there is no canonical normal vector field at the boundary defined from $g$ (recall that $g$ blows-up at 
the boundary), we can consider $\bar{g}:=\rho^2g$ for some geodesic boundary defining function and take the unit normal 
vector field for $\bar{g}$, that is $\nabla^{\bar{g}}\rho$, which we denote by $\pl_\rho$. 
It turns out that $(\pl_\rho^ku|_{\pl\bar{X}})_{k=1,\dots,n-1}$ are locally determined 
by $u|_{\pl\bar{X}}=f$ and the first term in the Taylor expansion of $u$ which is global is the $n$-th
$\pl_\rho^nu|_{\pl\bar{X}}.$ 
We thus define $\mc{S}f\in C^{\infty}(\pl\bar{X})$ by
\begin{equation}\label{defsc}
\mc{S}f:=\frac{1}{n!}\pl_\rho^nu|_{\pl\bar{X}}.
\end{equation}   
Notice that $\mc{S}$ a priori depends on the choice of $\rho$, we shall say that it is associated to $\rho$. 
It can be checked that if $\hat{\rho}=e^{\omega}\rho$ is another geodesic boundary 
defining function with $\omega\in C^{\infty}(\bar{X})$, then the scattering map 
$\hat{\mc{S}}$ associated to $\hat{\rho}$ satisfy $\hat{\mc{S}}=e^{-n\omega_0}\mc{S}$ where $\omega_0=
\omega|_{\pl\bar{X}}$, see \cite{GRZ} and Subsection \ref{scmap} below.

We also remark that the fact that $u\in C^\infty(\bar{X})$ strongly depends on the fact that the manifold under consideration is Einstein and has even dimensions.  For more general asymptotically hyperbolic manifolds, the solution $u$ to \eqref{poisson} possibly has a logarithmic singularity as shown in \cite{GRZ}.
Our second result is the following
\begin{theo}\label{th3}
Let $(X_1,g_1)$ and $(X_2,g_2)$ be connected, $C^\infty,$ $(n+1)$-dimensional conformally compact manifolds, with $n+1$ even. Suppose that $g_1$ and $g_2$ are Einstein and
that $\pl\bar{X_1}$ and $\pl\bar{X}_2$ contain a common open set $\Gamma$ such that the identity
map ${\rm Id}:\Gamma\subset\bar{X}_1\to \Gamma\subset \bar{X}_2$ is a smooth diffeomorphism.
If for $i=1,2$,
there exist boundary defining functions $\rho_i$ of $\partial \bar{X}_i$ such that 
the scattering maps $\mc{S}_i$ of $\Delta_{g_i}$ associated to 
$\rho_i$ satisfy $(\mc{S}_1f)|_{\Gamma}=(\mc{S}_2f)|_{\Gamma}$ for all $f\in C_0^\infty(\Gamma)$.
Then there is a diffeomorphism $J:\bar{X}_1\to\bar{X}_2$, such that  $J^{*}g_2=g_1$ in $X_1$.  
\end{theo}

The proofs  are based on the results of Lassas and Uhlmann \cite{LU},  and   Lassas, Taylor and Uhlmann \cite{LUT}, and  suitable unique continuation theorems for Einstein equation. 

It is shown in \cite{LU} that a connected compact manifold with boundary $(\bar{X}=X\cup\pl\bar{X},g)$, is determined by the Dirichlet-to-Neumann  if the interior $(X,g)$ is real analytic, and if there exists an open set $\Gamma$ of the boundary $\pl\bar{X}$ which is real analytic with $g$ real analytic up to $\Gamma$.  
In \cite{LUT} Lassas, Taylor and Uhlmann prove the analogue of this result for complete manifolds.

A theorem of De Turck and Kazdan, Theorem 5.2 of  \cite{DTK}, says that if $(\bar{X},g)$ is a connected Einstein manifold with boundary then the collection of harmonic coordinates give $X,$ the interior of $\bar{X},$  a real analytic structure which is compatible with its $C^\infty$ structure, and moreover $g$ is real analytic in those coordinates. The principle is that Einstein's equation becomes a non-linear elliptic system  with real 
analytic coefficient in these coordinates, thus the real analyticity of the metric. But since the harmonic coordinates satisfy the Laplace  equation, they are analytic as well.

However this construction is not necessarily valid at the boundary. Therefore one cannot guarantee that $(\bar{X},g)$ is real analytic at the boundary, and hence one cannot directly apply the results of \cite{LUT}. 

To prove Theorems \ref{th1} and \ref{th3}, we first 
show that the DN map (or the scattering map) determines the metric in a small 
neighbourhood $U$ of a point $p\in\Gamma\subset\pl\bar{X}$, then we shall prove that this 
determines the Green's function in $U\x U.$  However one of the  results of  
\cite{LUT} says that this determines the whole  Riemannian manifold, provided it is 
 real-analytic, but as mentioned above, this is the case of the interior of an Einstein manifold.

The essential part  in this paper is the reconstruction near the boundary.  This
will be done using the ellipticity of Einstein equation in harmonic coordinates,
and by applying a unique continuation theorem for the Cauchy problem for elliptic systems with diagonal principal part. The unique continuation result we need in the compact case was essentially proved by Calder\'on \cite{Ca,Ca1}.  The
conformally compact case is more involved since the system is  only elliptic in the
uniformly degenerate sense of \cite{Ma, Ma1, Ma2,MM}, see also \cite{Albao}.   
When the first version of this paper was completed we learned that  O. Biquard \cite{Bi} proved a unique continuation result for Einstein manifolds without using the DN map for functions, which was a problem that was part of the program of M. Anderson \cite{An}. 
Under our assumptions, it seems somehow natural to use harmonic coordinates for Einstein equation, and
we notice that our approach is self-contained and does not require the result of \cite{Bi}.

 Throughout this  paper when we refer to the real analyticity of the metric, we mean that it is real analytic with respect to the real analytic structure defined from harmonic coordinates corresponding to the Einstein metric $g.$

\section{Acknowledgments}
The work of both authors was funded by the NSF under grant DMS-0500788.
Guillarmou  acknowledges support of french ANR grants JC05-52556 and JC0546063 and thanks 
the MSI at ANU, Canberra, where part of this work was achieved. We thank
Erwann Delay, Robin Graham and Gunther Uhlmann for helpful conversations.  
Finally, we are grateful to the anonymous referee for a very careful reading, for many useful comments, and especially for uncovering an error in the first submitted version of the paper

\section{Inverse problem for Einstein manifolds with boundary}\label{inverseprob1}

The result of De Turck and Kazdan  concerning the analyticity of the metric does not apply to Einstein manifolds with boundary $(\bar{X}=X\cup\pl\bar{X},g).$  Their argument breaks down since the boundary can have 
low regularity even though $g$ has constant Ricci curvature. 
This means that the open incomplete manifold $(X,g)$ is real-analytic with respect to the analytic structure defined by harmonic coordinates,
but a priori $(\bar{X},g)$ does not satisfy this property.   We will use the Dirichlet-to-Neumann map to overcome this difficulty.

\subsection{The Dirichlet-to-Neumann map}\label{dnmap}
As in section \ref{intro}, $$\mc{N}: C^\infty(\pl\bar{X})\to C^{\infty}(\pl \bar{X})$$ 
is defined by solving the Dirichlet problem \eqref{solving} with $f\in C^{\infty}(\pl\bar{X}),$  and setting 
$\mc{N}f:=-\pl_{n}u|_M$ where $\pl_n$ is the interior pointing 
normal vector to the boundary for the metric $g$.
Its Schwartz kernel is related to the Green function $G(z,z')$ 
of the Laplacian $\Delta_g$ with Dirichlet condition on $\pl\bar{X}$
 by the following identity
\begin{lem}\label{kernel}
The Schwartz kernel $\mc{N}(y,y')$ of $\mc{N}$ is given for $y,y'\in \pl\bar{X}$, $y\not=y'$, by 
\[\mc{N}(y,y')=\pl_{n}\pl_{n'}G(z,z')|_{z=y,z'=y'}\]
where $\pl_{n},\pl_{n'}$ are respectively the inward pointing normal vector fields to the boundary in variable $z$ and $z'$. 
\end{lem}
\noindent \textsl{Proof}: Let $x$ be the distance function to the boundary in $\bar{X}$, it is smooth
in a neighbourhood of $\pl\bar{X}$ and the normal vector field to the boundary is the gradient 
$\pl_n=\nabla^g x$ of $x$. The flow $e^{t\pl_n}$ of $\nabla^g x$ induces a diffeomorphism
$\phi:[0,\eps)_t\x \pl\bar{X}\to \phi([0,\eps)\x\pl\bar{X})$ defined by $\phi(t,y):=e^{t\pl_n}(y)$ and we have $x(\phi(t,y))=t$. This induces natural coordinates $z=(x,y)$ near the boundary, these are normal geodesic coordinates. 
The function $u$ in \eqref{solving} can be obtained by taking 
\[u(z):=\chi(z) - \int_{\bar{X}} G(z,z')(\Delta_g \chi)(z')dz'\]
where $\chi$ is any smooth function on $\bar{X}$ such that $\chi=f+O(x^2)$. 
Now using Green's formula and $\Delta_g(z)G(z,z')=\delta(z-z')=\Delta_g(z')G(z,z')$
where $\delta(z-z')$ is the Dirac mass on the diagonal, we obtain for $z \in X$
\[\begin{split}
u(z)=&\int_{\pl\bar{X}}\Big(\pl_{n'}G(z,z')\chi(z')\Big)|_{z'=y'}dy'-\int_{\pl\bar{X}}
\Big(G(z,z')(\pl_n\chi)(z')\Big)|_{z'=y'}dy'\\
 u(z)=& \int_{\pl\bar{X}}\Big(\pl_{n'}G(z,z')\Big)|_{z'=y'}f(y')dy'\end{split}.
\]
We have a Taylor expansion $u(x,y)=f(y)+x\mc{N}f(y) +O(x^2)$ near the boundary.  Let $y\in \p X$ and take $f\in C^\infty(X)$ supported near $y.$ Thus pairing with $\phi\in C^{\infty}(\pl\bar{X})$ gives 
\begin{gather}
\int_{\pl\bar{X}}u(x,y)\phi(y)dy= \int_{\pl\bar{X}}f(y)\phi(y)dy-x\int_{\pl\bar{X}}\phi(y)\mc{N}f(y)dy+O(x^2).\label{neweq}
\end{gather}
Now taking $\phi$ with support disjoint to the support of $f,$ thus $\phi f=0,$ and differentiating  \eqref{neweq} in $x,$ we see, using the fact that Green's function $G(z,z')$ is smooth outside the 
diagonal, that 
\[\int_{\pl\bar{X}}\phi(y)\mc{N}f(y)dy=\int_{\pl\bar{X}}\int_{\pl\bar{X}}\Big(\pl_n\pl_{n'}G(z,z')\Big)|_{z=y,z'=y'}f(y')\phi(y)dydy',\] 
which proves the claim. 
\qed 

\subsection{The Ricci tensor in harmonic coordinates and unique continuation}\label{riccitensor}
Let us take coordinates $x=(x_0,x_1,\dots,x_n)$  near a point $p\in\pl \bar{X}$,
with $x_0$ a boundary defining function of $\pl\bar{X}$, then 
$\ric(g)$ is given by definition by 
\begin{equation}\label{defric}
\ric(g)_{ij}=\sum_{k}\Big( \pl_{x_k}\Gamma_{ji}^k-\pl_{x_j}\Gamma_{ki}^k+\sum_{l}\Gamma_{kl}^k\Gamma_{ji}^l
-\sum_l\Gamma_{jl}^k\Gamma_{ki}^l\Big)
\end{equation}
with 
\begin{equation}\label{christoff}
\Gamma_{ji}^k=\demi\sum_{m}g^{km}\Big( \pl_{x_i}g_{mj}+\pl_{x_j}g_{mi}-\pl_{x_m}g_{ij}\Big).  
\end{equation}
Lemma 1.1 of \cite{DTK} shows that $\Delta_gx_k=\sum_{i,j}g^{ij}\Gamma^k_{ij}$, so Einstein equation 
$\ric(g)=\la g$ for some $\la\in\rr$ can be written as the system (see also Lemma 4.1 in \cite{DTK})
\begin{equation}\label{ricci}
-\demi \sum_{\mu,\nu}g^{\mu\nu}\pl_{x_\mu}\pl_{x_\nu} g_{ij}+\demi\sum_{r}
(g_{ri}\pl_{x_j}(\Delta_gx_r) +g_{rj}\pl_{x_i}(\Delta_gx_r))+ Q_{ij}(g,\pl g)=0
\end{equation}
where $Q_{ij}(A,B)$ is smooth and polynomial of degree two in $B$, where $A,B$  
denote vectors $(g_{kl})_{k,l}\in \rr^{(n+1)^2}$ and $(\pl_{x_m}g_{kl})_{k,l,m}\in \rr^{(n+1)^3}.$
From this discussion   we deduce the following
\begin{prop}\label{system0}
Let $(x_0,x_1,\dots,x_n)$ be harmonic coordinates for $\Delta_g$ near a point $p\in\{x_0=0\}$, 
then there exist $Q_{ij}(A,B)$ smooth, polynomial of degree $2$ in $B\in\rr^{(n+1)^3}$, 
such that $\ric(g)=\la g$ is equivalent near $p$ to the system
\begin{equation}\label{ricci1}
\sum_{\mu,\nu}g^{\mu\nu}\pl_{x_\mu}\pl_{x_\nu} g_{ij}+ Q_{ij}(g,\pl g)=0,\quad i,j=0,\dots,n
\end{equation}
with $\pl g:=(\pl_{x_m}\bar{g}_{kl})_{k,l,m}\in \rr^{(n+1)^3}$.
\end{prop}
Now we may use a uniqueness theorem for the Cauchy problem of such elliptic systems.
\begin{prop}\label{calderon}
If $C:=(c_{ij})_{i,j=0,\dots,n},D:=(d_{ij})_{i,j=0,\dots,n}$ are smooth symmetric $2$-tensors near 
$p\in \{x_0=0\},$ with $C$ positive definite, the system \eqref{ricci1} near $p$ with boundary conditions $g_{ij}|_{x_0=0}=c_{ij}$ and $\pl_{x_0}g_{ij}|_{x_0=0}=d_{ij}$, $i,j=0,\dots,n$, has at most a unique smooth solution.
\end{prop}
\noindent \textsl{Proof}: The system is elliptic and the leading symbol is a scalar times the identity, the result could then be proved  using Carleman estimates. For instance, uniqueness properties are proved by Calderon \cite{Ca,Ca1} for elliptic systems when the characteristics of the system are non-multiple, but in our case they are multiple. However, since the leading symbol is scalar and this scalar symbol has only non-multiple characteristics, the technics used in Calderon could be applied like in the case of a single equation with non-multiple characteristics. 
Since we did not find references that we can cite directly, we prefer to 
use Proposition \ref{unique1} which is a consequence of 
a uniqueness result of Mazzeo \cite{Ma}. Indeed, first it is straightforward to notice, by using boundary normal coordinates,  that two solutions of \eqref{ricci1} with same Cauchy data
agree to infinite order at the boundary, therefore we may multiply \eqref{ricci1} by $x_0^2$ and \eqref{ricci1} becomes of the form \eqref{einsteing} thus Proposition \ref{unique1} below proves uniqueness.
\qed

\subsection{Reconstruction near the boundary}
Throughout this section we assume that $(\bar{X}_1,g_1)$, $(\bar{X}_2,g_2)$ are $C^\infty$ connected
Einstein manifolds with  boundary $M_j=\pl\bar{X}_j,$ $j=1,2,$ such that $M_1$ and $M_2$ contain a common open set $\Gamma,$ and that the identity map $\operatorname{Id} : \Gamma \subset \partial X_1 \longrightarrow
\Gamma \subset \partial X_2$ is a $C^\infty$ diffeomorphism.
Moreover we assume that for every $f\in C_0^\infty(\Gamma),$ the
 Dirichlet-to-Neumann maps satisfy $$\mc{N}_1f|_\Gamma=\mc{N}_2f|_\Gamma.$$ We first prove 
\begin{lem}\label{collar}
For $i=1,2$, there exists $p\in \Gamma$, some neighbouroods $U_i$ of $p$ in $\bar{X}_i$ and a diffeomorphism $F:U_1\to U_2$, $F|_{U_1\cap X_1}$ analytic, such that $F^*g_2=g_1$ and  $F|_{U_1\cap \Gamma}={\rm Id}$.
\end{lem}
\noindent \textsl{Proof}:  For $i=1,2$, let $t_i={\rm dist}(.,\pl\bar{X}_i)$ be the distance to the boundary in $\bar{X}_i$, then the flow $e^{t\nabla^{g_i}t_i}$ of the gradient $\nabla^{g_i} t_i$ induces a diffeomorphism
\begin{gather*}
\phi^i:[0,\eps)\x \pl\bar{X}_i \to \phi^i([0,\eps)\x \pl\bar{X}_i)\\
\phi^i(t,y):=e^{t\nabla^{g_i}t_i}(y),
\end{gather*}
and we have the decomposition near the boundary $(\phi^i)^* g_i=dt^2+h_i(t)$ for some one-parameter family of metrics $h_i(t)$ on $\pl\bar{X}_i$.  
Lee-Uhlmann \cite{LU} proved that $\mc{N}_1|_\Gamma=\mc{N}_2|_{\Gamma}$ implies that  
\begin{equation}\label{leeuhl}
\pl_{t}^kh_1|_\Gamma=\pl_{t}^kh_2|_\Gamma, \quad \forall \;\ k\in\nn_0.
\end{equation}
Let us now consider $H_i:={\phi^i}^*g_i$ on the collar $[0,\eps)_t\x \Gamma$.   
Let $p\in \Gamma$ be a point of the boundary and $(y_1,\dots,y_n)$ be a set of local coordinates in a neighbourhood of $p$
in $\Gamma,$  and extend each $y_j$ to $[0,\eps)\x \Gamma$ by the function $(t,m)\to y_j(m)$. 
Notice that $\phi^2\circ (\phi^1)^{-1}$ is a smooth diffeomorphism from a neighbourhood of $p$ in $\bar{X}_1$ to a neighbourhood of $p$ in $\bar{X}_2$, this is a consequence of the fact that ${\rm Id}:\Gamma\subset \bar{X}_1\to \Gamma\subset\bar{X}_2$ is a diffeomorphism.
Using $z:=(t,y_1,\dots, y_n)$ as coordinates on $[0,\eps)\x \Gamma$ near $p$, then \eqref{leeuhl} shows that there is  an open neighbourhood $U$ of $p$ in $[0,\eps)_t\x \Gamma$ such that $H_2=H_1+O(t^{\infty})$ and we can always
assume $U\cap \{t=0\}\not=\Gamma$. Let $y_0\in C_0^\infty(\Gamma)$ with $y_0=0$ on $U\cap\{t=0\}$ 
but $y_0$ not identically $0$, and by cutting off far from $p$ me may assume that $y_j\in C_0^\infty(\Gamma)$ for $j=1,\dots,n$.
Now let $(x_0^1,x^1_1,\dots,x^1_n)$ and $(x_0^2,x^2_1,\dots,x^2_n)$ be harmonic functions near $p$ in $[0,\eps)\x \Gamma$ 
for respectively $H_1$ and $H_2$ such that $x^1_j=x^2_j=y_j$ on $\{t=0\}$. 
These functions are constructed by solving the Dirichlet problem $\Delta_{g_i}w^i_j=0$ on $\bar{X}_i$ with boundary data
$w^i_j|_{M_i}=y_j,$ $i=1,2,$ and $j=0,\dots,n,$ and by setting $x^i_j={\phi^i}^*w^i_j$. Note that $\{m\in U\cap \Gamma; x_0^i(m)=0, d x_0^i(m)=0\}$
is a closed  set with empty interior in $U\cap\{t=0\},$ since otherwise $x_0^i$ would vanish to order $2$ on an open set of $\{t=0\}$, thus by unique continuation it would be identically $0$ since it is harmonic.
Then $(x_0^1,\dots,x_n^1)$ and $(x^2_0,\dots,x_n^2)$ form smooth coordinate systems near at least a common point 
of $U\cap\{t=0\}$; for convenience let us denote again $p$ this point and $U\subset [0,\eps)\x \Gamma$ an 
open set containing $p$ where they both form smooth coordinates. 

We have $\Delta_{H_1}(x_j^1-x_j^2)=O(t^\infty)$ and $\pl_{t}x_j^1|_{t=0}=\pl_{t}x_j^2|_{t=0}$ for all $j$ since $\mc{N}_1|_\Gamma=\mc{N}_2|_\Gamma$.
Since $u=x_j^1-x_j^2$ is solution of $\Delta_{H_1}u=O(t^\infty)$ in $U$ with $u$ vanishing at order $2$ at the boundary $t=0$, a standard Taylor expansion argument shows that $x_j^1=x_j^2+O(t^\infty)$ in $U$ for all $j$.
Now define $\psi:U\to \psi(U)\subset U$ 
so that $(x^1_0,\dots,x^1_n)=(\psi^*x_0^2,\dots,\psi^*x^2_n)$.
Then $\psi={\rm Id}+O(t^\infty)$ in $U$, and consequently we obtain in $U$ 
\begin{equation}\label{g2=g1}
\psi^*H_2=H_1+O(t^\infty).
\end{equation} 
The metrics $g=H_1$ and $g=\psi^*H_2$ both satisfy Einstein equation $\ric(g)=\la g$ in $U$. Moreover in coordinates $(x^1_0,\dots, x_n^1)$ this correspond to the system \eqref{ricci1} and since the coordinates are harmonic 
with respect to $g$, the system is elliptic and diagonal to leading order. 
From the unique continuation result in Proposition \ref{calderon},
we conclude that there exists a unique solution to this system in $U_1$ with given initial data $g|_{U\cap \{t=0\}}$ and $\pl_{x^1_0}g|_{U\cap\{t=0\}}$. In view of \eqref{g2=g1}, this proves that $H_1=\psi^*H_2$ in $U$. 
Although it is not relevant for the proof, we remark that $\psi$ is actually the Identity on $U$ since $\psi|_{U\cap\Gamma}={\rm Id}$ and 
it pulls back one metric in geodesic normal coordinates to the other.  
Now it suffices to go back to $\bar{X}_1$ and $\bar{X}_2$ through $\phi^1,\phi^2$ and we have proved the Lemma
by setting $U_i:=\phi^i( U)$ and 
\begin{gather}
F:=\phi^2\circ\psi\circ(\phi^1)^{-1}.  \label{defF}
\end{gather}
Remark that $F$ is analytic from $U_1\cap \{t_1\not=0\}$ to $U_2\cap\{t_2\not=0\}$ 
since the harmonic functions $w_j^i$ define the analytic structure in $U_i\cap\{t_i\not=0\}$ for all $j=0,\dots,n$ and 
$F$ is the map that identify $w^1_j$ to $w_j^2$ for all $j$.\qed\\

Next we prove 
\begin{cor}\label{corol1}
For $i=1,2$, let $G_i(z,z')$ be the Green function of  $\Delta_{g_i}$ in
$\bar{X}_i$ with Dirichlet condition at $M$, then $\mc{N}_1|_\Gamma=\mc{N}_2|_\Gamma$ implies that there exists
an open set $U_1\subset X_1$ with
\[G_2(F(z),F(z'))=G_1(z,z'), \quad (z,z')\in (U_1\x U_1)\setminus \{z=z'\},\] 
where $F$ was defined in \eqref{defF}
\end{cor}
\noindent \textsl{Proof}: First we remark that $g_1$ is Einstein and thus real analytic in $U_1\setminus (U_1\cap M)$, 
so is any harmonic function in this open set. 
Let $\pl_n,\pl_{n'}$ be the normal vector fields to the boundary in the first and second variables in $U_1\x U_1$ respectively, as
defined in Lemma \ref{kernel}. We see from the proof of Lemma \ref{collar} that $F_*\pl_n$ and 
$F_*\pl_{n'}$ are the normal vector fields to the boundary in the first and second variable in $U_2\x U_2$ 
(since $\psi={\rm Id}+O(t^\infty)$ in that Lemma). So we get $\pl_{n'}G_2(F(z),F(z'))=(F_*\pl_{n'})
G_2(F(z),z')$ for $z'\in M$ since $F|_{U_1\cap\Gamma}={\rm Id}$.

We first show that $\pl_{n'}G_2(F(z),F(z'))=\pl_{n'}G_1(z,z')$
for any $(z,z')\in U_1\x (U_1\cap \Gamma)\setminus\{z=z'\}$.
Now fix $z'\in U_1\cap \Gamma$, then the function $T_1(z):=\pl_{n'}G_1(z,z')$ solves $\Delta_{g_1}T_1=0$ in $
U_1\setminus\{z'\}$ and, using Lemma \ref{kernel}, it has boundary values $T_1|_{U_1\cap \Gamma\setminus\{z'\}}=0$ 
and $\pl_nT_1|_{U_1\cap \Gamma\setminus\{z'\}}=\mc{N}_1(.,z')$ where $\mc{N}_i(.,.)$ denote the Schwartz kernel of 
$\mc{N}_i$, $i=1,2$.
The function $T_2(z):=\pl_{n'}G_{2}(F(z),F(z'))$ solves
$\Delta_{F^*g_2}T_2(z)=\Delta_{g_1}T_2(z)=0$ in $U_1\setminus\{z'\}$.
We also have $\pl_{n}T_2|_{U_1\cap \Gamma\setminus\{z'\}}=
F^*[(F_*\pl_{n})(F_*\pl_{n'})G_2(.,z')|_{U_1\cap \Gamma\setminus\{z'\}}]$ and $T_2|_{U_1\cap \Gamma\setminus\{z'\}}=0$. 
But from Lemma \ref{kernel},  
$(F_*\pl_n)(F_*\pl_{n'})G_2(.,z')|_{U_1\cap \Gamma\setminus\{z'\}}=\mc{N}_2(.,z')$ where $\mc{N}_2(.,.)$ 
is the Schwartz kernel of $\mc{N}_2$. Using again that $F|_{U_1\cap\Gamma}={\rm Id}$, we deduce that
$\pl_{n}T_2|_{U_1\cap \Gamma\setminus\{z'\}}=\mc{N}_2(.,z')$. 
By our assumption $\mc{N}_1|_\Gamma=\mc{N}_2|_\Gamma$, we conclude that $T_1$ and $T_2$ solve the same Cauchy problem  
near $U_1\cap \Gamma\setminus\{z'\}$, so first by unique continuation near the boundary and 
then real analyticity in $U_1\setminus(U_1\cap \Gamma)$, we obtain $T_1=T_2$ there. 

Now we can use again similar arguments to prove that $G_1(z,z')=G_2(F(z),F(z'))$ in $(U_1\x U_1)\setminus
\{z=z'\}$. Indeed, fix $z'\in U_1$, then $T_1(z):=G_1(z',z)$ and $T_2(z):=G_2(F(z'),F(z))$ solve
$\Delta_{g_1}T_i=0$ in $U_1\setminus\{z'\}$ and with boundary values $T_i|_\Gamma=0$ and $\pl_nT_1|_{U_1\cap\Gamma}=\pl_nT_2|_{U_1\cap\Gamma}$
by what we proved above. Thus unique continuation for Cauchy problem and real analyticity allow us to conclude that
$T_1=T_2$.
\qed

\subsection{Proof of Theorem \ref{th1}}\label{proofth1}
To conclude the proof of \ref{th1}, we use the following Proposition which is implicitly
proved by Lassas-Taylor-Uhlmann \cite{LUT}
\begin{prop}\label{lut1}
For $i=1,2$, let $(\bar{X}_i,g_i)$ be $C^\infty$ connected Riemannian manifolds with boundary, 
assume that its interior
$X_i$ has a real-analytic structure compatible with the smooth structure and such that the
 metric $g_i$ is real analytic on $X_i$. Let $G_i(z,z')$ be the Green function of the Laplacian
$\Delta_{g_i}$ with Dirichlet condition at $\pl\bar{X}_i$, and assume there exists an open set 
$U_1\subset X_1$ and an analytic diffeomorphism $F:U_1\to F(U_1)\subset X_2$ such that $G_1(z,z')=G_2(F(z),F(z'))$
for $(z,z')\in (U_1\x U_1)\setminus\{z=z'\}$. Then there exists a diffeomorphism $J:X_1\to X_2$ such that
$J^*g_2=g_1$ and $J|_{U_1}=F$. 
\end{prop}

The proof is entirely done in section 3 of \cite{LUT}, although not explicitly written under that form.
The principle is to define maps 
\[\mc{G}_j: X_j\to H^s(U_1), \quad \mc{G}_1(z):=G_1(z,.), \quad \mc{G}_2(z):=G_2(z,F(.))\]
where $H^s(U_1)$ is the $s$-Sobolev space of $U_1$ for some $s<1-(n+1)/2$, then prove 
that $\mc{G}_j$ are embeddings with $\mc{G}_1(X_1)=\mc{G}_2(X_2)$, 
and finally show that $J:=\mc{G}_2^{-1}\circ\mc{G}_1:X_1\to X_2$ is an isometry.
Note that $J$ restricts to $F$ on $U_1$ since $G_1(z,z')=G_2(F(z),F(z'))$.\\ 

Proposition \ref{lut1} and Corollary \ref{corol1} imply Theorem \ref{th1}, after noticing that an isometry 
$\psi:(X_1,g_1)\to (X_2,g_2)$ extends smoothly to the manifold with boundary $(\bar{X}_1,g_1)$
by smoothness of the metrics $g_i$ up to the boundaries $\pl\bar{X}_i$.

\section{Inverse scattering for conformally compact Einstein manifolds}

Consider an $n+1$ dimensional connected conformally compact Einstein manifold $(\bar{X},g)$ with $n+1$ even, 
and let $\rho$ be a geodesic boundary defining function and $\bar{g}:=\rho^2g$. 
Using the flow $\phi_t(y)$ of the gradient $\nabla^{\rho^2g}\rho$, one has a diffeomorphism $\phi:[0,\eps)_t\x \pl\bar{X}
\to \phi([0,\eps)\x \pl\bar{X})\subset\bar{X}$ defined by $\phi(t,y):=\phi_t(y)$, and the metric pulls back to
\begin{equation}\label{model}
\phi^*g=\frac{dt^2+h(t)}{t^2}
\end{equation}
for some smooth one-parameter family of metrics $h(t)$ on the boundary $\pl\bar{X}$. Note that here $\phi^*\rho=t$.

\subsection{The scattering map}\label{scmap}
The scattering map $\mc{S},$ or scattering matrix, defined in the introduction is really $\mc{S}=S(n)$, where $S(\la)$ for $\la\in\cc$
is defined in \cite{JSB,GRZ}.  
Let us construct $\mc{S}$ by solving the boundary value problem $\Delta_gu=0$ with $u\in C^{\infty}(\bar{X})$
and $u|_{\pl\bar{X}}=f$ where $f\in C^{\infty}(\pl\bar{X})$ is given.   This follows the construction in section 4.1 of \cite{GRZ}.
Writing $\Delta_g$ in the collar $[0,\eps)_t\x\pl\bar{X}$ through the diffeomorphism $\phi$, we have
\[\Delta_g=-(t\pl_t)^2+(n-\frac{t}{2}\tra_{h(t)}(\pl_th(t)))
t\pl_t+t^2\Delta_{h(t)}\]
and for any $f_j\in C^\infty(\pl\bar{X})$ and $j\in\nn_0$
\begin{equation}\label{indiciel}
\Delta_g(f_j(y)t^j)=j(n-j)f_j(y)t^{j}+t^{j}(H(n-j)f_j)(t,y),
\end{equation}
\[(H(z)f_j)(t,y):=t^2\Delta_{h(t)}f_j(y)-
\frac{(n-z)t}{2}\tra_{h(t)}(\pl_th(t))f_j(y).\]
Now recall that since $g$ is Einstein and even dimensional, we have $\pl^{2j+1}_th(0)=0$ for $j\in\nn_0$ such that
$2j+1<n$, see for instance Section 2 of \cite{GR}. Consequently, $H(n-j)f_j$ is an even function of $t$ modulo
$O(t^n)$ for $j\not=0$, and modulo $O(t^{n+2})$ when $j=0$.
Since $H(n-j)f_j$ also vanishes at $t=0$, we can construct by induction a Taylor series using \eqref{indiciel} 
\begin{equation}\label{constfj}
F_j=\sum_{k=0}^jt^{k}f_{k}(y), \quad F_0=f_0=f, \quad F_j=F_{j-1}+t^j\frac{[t^{-j}(\Delta_gF_{j-1})]|_{t=0}}{j(j-n)} 
\end{equation}
for $j<n$ such that $\Delta_gF_j=O(t^{j+1})$. Note that, since $H(n-j)f_j$ has even powers of $t$ modulo $O(t^n)$,
we get $f_{2j+1}=0$ for $2j+1<n$.
For $j=n$, the construction of $F_n$ seems to fail but actually we can remark
that $\Delta_{g}F_{n-1}=O(t^{n+1})$ instead of $O(t^n)$ thanks to the fact that 
$t^{2j}H(n-{2j})f_{2j}$ has even Taylor expansion at $t=0$ modulo $O(t^{2j+n+2})$ 
by the discussion above. So we can set $F_n:=F_{n-1}$ and then continue to define $F_j$ for $j>n$ using 
\eqref{constfj}. Using Borel's Lemma, one can construct $F_{\infty}\in C^{\infty}(\bar{X})$ such that
$\phi^{*}F_{\infty}-F_j=O(t^{j+1})$ for all $j\in\nn$ and $\Delta_gF_{\infty}=O(\rho^{\infty})$.
Now we finally set $u=F_{\infty}-G\Delta_gF_\infty$ where $G:L^2(X,{\rm dvol}_g)\to L^2(X,{\rm dvol}_g)$ 
is the Green operator, i.e. such that $\Delta_gG={\rm Id}$, recalling that $\ker_{L^2}\Delta_g=0$ by \cite{MM}. 
From the analysis of $G$ in \cite{MM}, one has that $G$ maps $\dot{C}^{\infty}(\bar{X}):=\{v\in 
C^{\infty}(\bar{X}), v=O(\rho^{\infty})\}$ to $\rho^{n}C^{\infty}(\bar{X})$. This proves that 
$u\in C^{\infty}(\bar{X})$ and has an asymptotic expansion
\begin{equation}\label{asympu}
\phi^*u(t,y)=f(y)+\sum_{0<2j<n}t^{2j}f_{2j}(y)-\phi^*(G\Delta_gF_{\infty})+O(t^{n+1}).
\end{equation}
In particular the first odd power is of order $t^n$ and its coefficient is given by
the smooth function $[t^{-n}\phi^*(G\Delta_gF_{\infty})]_{t=0}$ of $C^{\infty}(\pl\bar{X})$.
Notice that the $f_{2j}$ in the construction are local with respect to $f$, more precisely
$f_{2j}=p_{2j}f$ for some differential operator $p_{2j}$ on the boundary.
Note that we used strongly that the Taylor expansion of the metric $t^2\phi^*g$ at $t=0$ is even to order $t^n$, which comes from the fact that $X$ is Einstein and has even dimensions. Indeed for a general asymptotically hyperbolic manifold, $u$ has logarithmic singularities, see \cite{GRZ,Gu}.

Since $\phi^*\nabla^{\rho^2g}\rho=\pl_t$, the definition of $\mc{S}f$ in 
the Introduction is equivalent to $\mc{S}f=\frac{1}{n!}\pl^{n}_t\phi^*u|_{t=0}$, i.e. the $n$-th Taylor 
coefficient of the expansion of $\phi^*u$ at $t=0$, in other words 
\[\mc{S}f=-[t^{-n}\phi^*(G\Delta_gF_{\infty})]_{t=0}=-[\rho^{-n}G\Delta_gF_\infty]|_{\pl\bar{X}}.\]  
From the analysis of Mazzeo-Melrose \cite{MM}, one can describe the 
behaviour of the Green kernel $G(z,z')$ near the boundary and outside 
the diagonal ${\rm diag}_{\bar{X}\x\bar{X}}$:
\begin{equation}\label{mazzeo}
\rho(z)^{-n}\rho(z')^{-n}G(z,z') \in C^{\infty}(\bar{X}\x\bar{X}\setminus{\rm diag}_{\bar{X}\x\bar{X}}).
\end{equation}
We can show easily that the kernel of $\mc{S}$ is the boundary value of \eqref{mazzeo} 
at the corner $\pl\bar{X}\x\pl\bar{X}$:
\begin{lem}\label{noyausc}
The Schwartz kernel $\mc{S}(y,y')$ of the scattering map $\mc{S}$ is, for $y\not=y'$,
\[\mc{S}(y,y')=n[\rho(z)^{-n}\rho(z')^{-n}G(z,z')]_{z=y,z'=y'}\]
where $G(z,z')$ is the Green kernel for $\Delta_g$.
\end{lem}
\noindent \textsl{Proof}:  Consider $(G\Delta_gF_\infty)(z)$ for $z\in X$ fixed and use Green formula
on the compact $U_\eps:=\{z'\in X; \rho(z)\geq \eps, {\rm dist}(z',z)\geq \eps\}$
\[\int_{U_\eps}G(z,z')\Delta_gF_{\infty}(z'){\rm dv}_{g}(z')=\int_{\pl U_\eps}
(G(z,z')\pl_{n'}F_{\infty}(z')-\pl_{n'}G(z,z')F_{\infty}(z'))d\nu_{\eps}(z')\]
where $\pl_{n'}$ is the unit normal interior pointing vector field of $\pl U_{\eps}$ (in the right 
variable $z'$) and 
$d\nu_{\eps}$ the measure induced by $g$ there. Consider   
the part $\rho(z')=\eps$ in the variables as in \eqref{model} using 
the diffeomorphism $\phi$, i.e. $\phi(t',y')=z'$, 
then $\phi^*\pl_{n'}=t'\pl_t'$ and $\phi^*(d\nu_{t'})={t}'^{-n}{\rm dvol}_{h(t')}$.
Using \eqref{mazzeo} and $F_{\infty}=f+O(\rho^2)$ by the construction of $F_{\infty}$ above 
the Lemma, we see that the integral on $\rho'=\eps$ converges to
\[n\int_{\pl{\bar{X}}}[\rho(z')^{-n}G(z,z')]_{z'=y'}f(y'){\rm dv}_{h(0)}(y').
\]
as $\eps\to 0.$ As for the part on ${\rm dist}(z',z)=\eps$, by another application of Green 
formula and $\Delta_g(z')G(z,z')=\delta(z-z')$, this converges to $F_\infty(z)$ as $\eps\to 0$.
We deduce that the solution $u$ of $\Delta_gu$ with $u|_{\pl\bar{X}}=f$ is given by 
\begin{equation}\label{uz}
u(z)=n\int_{\pl\bar{X}}[\rho(z')^{-n}G(z,z')]_{z'=y'}f(y'){\rm dv}_{h(0)}(y').
\end{equation}
Let us write $dy$ for ${\rm dv}_{h_0}(y)$. So given $y\in \p X,$ let $f$ be supported in a neighborhood of $y$ and take $\psi\in C^{\infty}(\pl\bar{X})$ with $\psi f=0$ and consider the pairing
\[\int_{\pl\bar{X}}\phi^{*}u(t,y)\psi(y)dy.\]
The Taylor expansion of $u$ at $t=0$ and the structure of $G(z,z')$ given by \eqref{mazzeo} show that 
\[\int_{\pl\bar{X}}\psi(y)\mc{S}f(y)dy=n\int_{\pl\bar{X}}[\rho(z)^{-n}\rho(z')^{-n}G(z,z')]|_{z=y,z'=y'}
\psi(y)f(y')dy'dy,\]
which proves the claim.
\qed\\

\textsl{Remark}: A more general relation between the kernel of the resolvent of $\Delta_g$, $(\Delta_g-\la(n-\la))^{-1}$, 
and the kernel of the scattering operator $S(\la)$ holds, as proved in \cite{JSB}.  But since 
the proof of Lemma \ref{noyausc} is rather elementary, we included it to make the paper
essentially self-contained.
 
\subsection{Einstein equation for $g$}\label{einsforg}
We shall analyze Einstein equation in a good system of coordinates, actually
constructed from harmonic coordinates for $\Delta_g$.
First choose coordinates $(y_1,\dots,y_n)$ in a neighbourhood $V\subset\pl\bar{X}$ 
of $p\in\pl\bar{X}$. Take an open set $W\subset \pl\bar{X}$ which contains $V$,
we may assume that $y_i\in C_0^\infty(W)$. 
Let $\phi$ be the diffeomorpism as in \eqref{model}. From the properties of the solution of the equation
$\Delta_gu=0,$ as in subsection \ref{scmap} (which follows Graham-Zworski \cite{GRZ}),  there exists $n$ smooth functions $(x_1,\dots,x_n)$ on $\bar{X}$ such that
\[\Delta_gx_i=0, \quad \phi^*x_i=y_i+\sum_{0<2k<n}t^{2k}p_{2k}y_i +t^n \mc{S}y_i+O(t^{n+1})
\]
where $p_k$ are differential operators on $\pl\bar{X}$ determined by the
 $(\pl^k_th(0))_{k=0,\dots,n-1}$ at the boundary (using the form \eqref{model}).
Similarly let $y_0\in C_0^\infty(W)$ be a non zero smooth function such that $y_0=0$ in $V$, 
then by Subsection \ref{scmap}  there exists $v \in C^{\infty}(\bar{X})$ such that
\[\Delta_g v=0, \quad \phi^*v=y_0+\sum_{0<2k<n}t^{2k}p_{2k}y_0 +t^n \mc{S}y_0+O(t^{n+1}).\]
Thus in particular $v$ vanishes near $p$ to order $\rho^n$ since $p_ky_0=0$ in $V$ for $k=1,\dots,n$, 
thus one can write
\[v=\rho^n(w +O(\rho)), \]
where $w$ is a smooth function on $\pl\bar{X}$ near $p$. The set $\{m \in V; w(m)\not=0\}$ is an open  
dense set of $V$. Indeed, otherwise $w$ would vanish in an open set of $V$ but an easy computation  
shows that if $U\in \rho^jC^{\infty}(\bar{X})$ then $\Delta_gU=-j(j-n)U+O(\rho^{j+1})$
so $v$ would vanish to infinite order at an open set of $V$ and by Mazzeo's unique continuation result \cite{Ma},
it would vanish identically in $\bar{X}$. 
Thus, possibly by changing $p$ to another point (still denoted $p$ for convenience), 
there exists $v\in C^{\infty}(\bar{X})$ such that $v$ is harmonic for $\Delta_g$ and 
$v=\rho^n(w+O(\rho))$ with $w>0$ near $p$, the function
$x_0:=v^{1/n}$ then defines a boundary defining function of $\pl\bar{X}$ near $p$, 
it can be written as $\rho e^f$ for some smooth $f$. 
Then $(x_0,x_1,\dots,x_n)$ defines a system of coordinates near $p$.\\
 
Let us now consider Einstein equations in these coordinates. Again like \eqref{ricci}, the principal part 
of $\ric(g)$ is given by
\[-\demi \sum_{\mu,\nu}g^{\mu\nu}\pl_{x_\mu}\pl_{x_\nu} g_{ij}+\demi\sum_{r}
(g_{ri}\pl_{x_j}(\Delta_gx_r) +g_{rj}\pl_{x_i}(\Delta_gx_r)).\]  
But all functions $x_r$ are harmonic, except $x_0,$ and the latter satisfies
\[0=\Delta_{g}x_0^n=-n\textrm{div}_g(x_0^{n-1}\nabla^g x_0)=-nx_0^{n-1}\Delta_{g}x_0-n(n-1)x_0^{n-2}|dx_0|_g^2\]
or equivalently $\Delta_gx_0=(1-n)x_0|dx_0|_{x_0^2g}$. But this involves only terms of order $0$ in the metric 
$g$ or $\bar{g}:=x_0^2g$ so the principal part of $\ric(g)$ in these coordinates is  
\[-\demi \sum_{\mu,\nu}g^{\mu\nu}\pl_{x_\mu}\pl_{x_\nu} g_{ij}\]
which is elliptic in the interior $X$.
We multiply the equation $\ric(g)=-ng$ by $x_0^2$ near $p$ and using (\ref{defric}) and (\ref{christoff}), with the 
commutations relations $[x_0\pl_{x_0},x_0^\alpha]=\alpha x_0^\alpha$ for all $\alpha\in\cc$, 
it is straightforward to obtain 
\begin{lem}\label{system2}
Let $x=(x_0,x_1,\dots,x_n)$ be the coordinates defined above near a point $p\in\{x_0=0\}$, 
then Einstein equation for $g$ can be written under the system
\begin{equation}\label{einsteing}
\sum_{\mu,\nu}x_0^2\bar{g}^{\mu\nu}\pl_{x_\mu}\pl_{x_\nu} \bar{g}_{ij}+ Q_{ij}(x_0,\bar{g},x_0\pl \bar{g})=0, \quad i,j=0,\dots,n
\end{equation}
where $\bar{g}=x_0^2g$ near $p$, $Q_{ij}(x_0,A,B)$ are smooth and polynomial of order $2$ in $B$, and $x_0\pl \bar{g}:=
(x_0\pl_{x_m} \bar{g}_{ij})_{m,i,j}\in\rr^{(n+1)^3}$.
\end{lem}

This is a non-linear system of order $2$, elliptic in the uniformly degenerate sense of \cite{Ma1,Ma2,MM} and 
diagonal at leading order.
We state the following unique continuation result for this system:
\begin{prop}\label{unique1}
Assume $\bar{g}_1$ and $\bar{g}_2$ are two smooth solutions of the 
system \eqref{einsteing} with $\bar{g}_1=\bar{g}_2+O(x_0^{\infty})$ near $p$. Then $\bar{g}_1=\bar{g}_2$
near $p$. 
\end{prop}
\noindent \textsl{Proof}: This is an  application of Mazzeo's unique continuation result \cite{Ma}.
We work in a small neighbourhood $U$ of $p$ and set $w=(\bar{g}_1-\bar{g}_2)$ near $p$. For $h$ metric near $p$ \
and $\ell$ symmetric tensor near $p$, let 
\[G(x_0,h,\ell):=-\sum_{\mu,\nu}x_0^2h^{\mu\nu}\pl_{x_\mu}\pl_{x_\nu}\bar{g}_{2}-
Q(x_0,h,\ell)\]
where $Q:=(Q_{ij})_{i,j=0,\dots,n}$. Note that $G$ is smooth in all its components. We have from \eqref{einsteing}
\begin{equation}\label{eqr2}
\sum_{\mu,\nu}x_0^2\bar{g}_1^{\mu\nu}\pl_{x_\mu}\pl_{x_\nu}w =G(x_0,\bar{g}_1,x_0\pl\bar{g}_1)-G(x_0,\bar{g}_2,x_0\pl\bar{g}_2).
\end{equation} 
Let $g_1:=x_0^{-2}\bar{g}_1$ and let $\nabla$ be the connection 
on symmetric $2$ tensors on $U$ induced  by $g_1$, then $\nabla^*\nabla w=\sum_{\mu,\nu}g_1^{\mu\nu}\nabla_{\pl_{x_\nu}}\nabla_{\pl_{x_\mu}}w$ and
in coordinates it is easy to check that 
$x_0(\nabla_{\pl_{x_\mu}}-\pl_{x_\mu})$ is a zeroth order operator with smooth coefficients 
up to the boundary, using \eqref{christoff} for instance. 
Therefore one obtains, using \eqref{eqr2}, 
\[|\nabla^*\nabla w|_{g_1}\leq C(|w|_{g_1}+|\nabla w|_{g_1})\]
for some $C$ depending on $\bar{g}_1, \bar{g}_2$. It then suffices to apply Corollary 11 of \cite{Ma}, this proves that
$w=0$ and we are done.
\qed

\subsection{Reconstruction near the boundary and proof of Theorem \ref{th3}}
The proof of Theorem \ref{th3} is  fairly close to that of Theorem \ref{th1}.
Let $(\bar{X}_1,g_1)$ and $(\bar{X}_2,g_2)$ be conformally compact Einstein manifolds with 
geodesic boundary defining functions $\rho_1$ and $\rho_2$. 
Let $\mc{S}_i$ be the scattering map for $g_i$ defined by \eqref{defsc} using the boundary defining function
$\rho_i$, assume that $\pl\bar{X}_1$ and $\pl\bar{X}_2$ contain a common open set $\Gamma$ such that the identity map which identifies the copies of $\Gamma$ is a diffeomorphism, and that
$\mc{S}_1f|_\Gamma=\mc{S}_2f|_{\Gamma}$ for all $f\in C^{\infty}_0(\Gamma)$.
Using the geodesic boundary defining function $\rho_i$ for $g_i$,$i=1,2$, there is a diffeomorphism 
$\phi^i:[0,\eps)_t\x\pl\bar{X}_i\to \phi^i([0,\eps)\x\pl\bar{X}_i)\subset \bar{X}_i$ as in \eqref{model}
so that 
\begin{equation}\label{modeli}
(\phi^i)^*g_i=\frac{dt^2+h_i(t)}{t^2}
\end{equation}
where $h_i(t)$ is a family of metric on $\pl\bar{X}_i$. We first show the 
\begin{lem}\label{infiniteorder}
The metrics $h_1(t)$ and $h_2(t)$ satisfy $\pl^j_th_1(0)|_{\Gamma}=\pl_t^jh_2(0)|_\Gamma$ for all $j\in\nn_0$.
\end{lem}
\noindent \textsl{Proof}: For a compact manifold $M,$ let us denote $\Psi^z(M)$ the set of classical pseudo-differential operators
of order $z\in\rr$ on $M$. 
Since $\mc{S}_i$ is the scattering operator $S_i(\la)$ at energy $\la=n$ for $\Delta_{g_i}$ as defined in \cite{JSB}, 
we can use \cite[Th.1.1]{JSB}, then we have that $\mc{S}_i\in\Psi^n(\pl\bar{X}_i)$ for $i=1,2$, with principal symbol 
$\sigma_n^i(y,\xi)=2^{-n}\Gamma(-\ndemi)/\Gamma(\ndemi)|\xi|_{h_i(0)}$, thus $h_1(0)=h_2(0)$ on 
$\Gamma$ and $\chi(\mc{S}_1-\mc{S}_2)\chi\in\Psi^{n+1}(\Gamma)$ for all $\chi\in C_0^\infty(\Gamma)$. 
Now we use Einstein equation, for instance the results of \cite{FGR,FGR2} (see also \cite[Sec. 2]{GR}) 
show, using only Taylor expansion of $\ric(g)=-ng$ at the boundary, that 
\[\pl^j_th_1(0)|_{\Gamma}=\pl^j_th_2(0)|_{\Gamma}, \quad j=0,\dots,n-1.\]
Then we use Theorem 1.2 of \cite{JSB} which computes the principal symbol of $\mc{S}_1-\mc{S}_2$. Since this result is entirely local, we can rephrase it on the piece $\Gamma$ of the boundary:  
if there exists a symmetric $2$-tensor $L$ on $\Gamma$ such that 
$h_1(t)=h_2(t)+t^kL+O(t^{k+1})$ on $[0,\eps)_t\x\Gamma$ 
for some $k\in\nn$, then for any $\chi\in C_0^\infty(\Gamma)$ we have 
$\chi(\mc{S}_1-\mc{S}_2)\chi\in\Psi^{n-k}(\Gamma)$ and the principal symbol of this operator at $(y,\xi)\in T^*\Gamma$ is
\footnote{It is important to notice that the coefficient of $|\xi|^{n-k}$ in \eqref{prsymb} is not exactly
that of Theorem 1.2 of \cite{JSB}, indeed there is a typo in equation (3.5) in \cite[Prop 3.1]{JSB}:
the coefficient in front of $T=\tra_{h_0}(L)$ there should be $k(k-n)/4$ instead of $k(k+1)/4$, this comes from the fact
that, in the proof of \cite[Prop 3.1]{JSB}, 
the term 
\[\frac{1}{16}x^2f\pl_x\log(\delta_2/\delta_1)\pl_x\log(\delta_2\delta_1)=-\frac{k(n+1)}{4}fx^{k}T+O(x^{k+1})\]
while it has been considered as a $O(x^{k+1})$ there.}
\begin{equation}\label{prsymb}
A_1(k,n)L(\xi^*,\xi^*)|\xi|^{n-k-2}_{h_0}+A_2(k,n)\frac{k(k-n)\tra_{h_0}(L)}{4}|\xi|_{h_0}^{n-k}
\end{equation}
where $h_0:=h_1(0)|_{\Gamma}=h_2(0)|_{\Gamma}$, $\xi^*:=h_0^{-1}\xi\in T_y\Gamma$ is the dual of 
$\xi$ through $h_0$, 
and $A_i(k,\la)$ are the meromorphic functions of $\la\in\cc$ defined by
\[A_1(k,\la):=-\pi^{-\ndemi}2^{k-2\la+n}\frac{\Gamma(\ndemi-\la+\frac{k}{2}+1)}{\Gamma(\la-\frac{k}{2}-1)}\frac{\Gamma(\la)^2}{\Gamma(\la-\ndemi+1)^2}\frac{T_1(k,\la)}{M(\la)}\]
\[A_2(k,\la):=\pi^{-\ndemi}2^{k-2\la+n-2}\frac{\Gamma(\ndemi-\la+\frac{k}{2})}{\Gamma(\la-\frac{k}{2})}\frac{\Gamma(\la)^2}{\Gamma(\la-\ndemi+1)^2}\frac{T_2(k,\la)}{M(\la)}\]
where $T_l(k,\la)$ is defined, when the integral converges, by  
\[T_l(k,\la):=\int_{0}^\infty\int_{\rr^n}\frac{u^{2\la-n+k+3-2l}}{(u^2+|v|^2)^{\la}(u^2+|e_1-v|^2)^\la}dv_{\rr^n}du, \quad e_1=(1,0,\dots,0),\]
and $M(\la)\in\cc$ is a constant not explicitly computed in \cite{JSB}. However at $\la=n$ the constant $M(n)$ is defined in \cite[Sec. 4]{JSB} such that $u(z)/n=M(n)f+O(\rho(z))$ where $u(z)$ is the function of \eqref{uz}, so $M(n)=n$ by \eqref{asympu}. 
Since we are interested in the case $k=n$, only the term with $A_1(n,n)$ appears, 
and setting $\la=n$ in $A_1(n,\la)$, with the explicit formulae above and the fact 
that $T_1(n,n)>0$ converges by Lemma 5.2 of \cite{JSB}, we see easily that $A_1(n,n)\not=0$ if $n>2$.   
Since we assumed $\chi \mc{S}_1\chi=\chi \mc{S}_2\chi$, this implies that $L=0$ and $h_1-h_2=O(t^{n+1})$ near $\Gamma$.
We finally use again \cite{FGR} (see \cite[Sect. 4]{FGR2} for full proofs), where it is proved that if $g_1=g_2+O(\rho^{n-1})$ with
$g_1,g_2$ conformally compact Einstein and $n$ odd, then $g_1=g_2+O(\rho^{\infty})$. Notice that their 
arguments are entirely local near any point of the boundary, so we can apply it near the piece 
$\Gamma$ of the boundary.
\qed\\

\begin{lem}\label{collar3}
For $i=1,2$, there exist $p\in\Gamma$,  neighbourhoods $U_i$ of $p$ in $\bar{X}_i$ and a diffeomorphism $F:U_1\to U_2$, $F|_{U_1\cap X_1}$ analytic, such that $F^*g_2=g_1$ and $F|_{U_1\cap\Gamma}={\rm Id}$. 
\end{lem}
\noindent \textsl{Proof}:  We  work in the collar $[0,\eps)_t\x \Gamma$ through the diffeomorphisms $\phi^i$ as in \eqref{modeli}. 
In a neighbourhood $U\subset [0,\eps)\x\Gamma$ of $p\in\Gamma$, we use coordinates $\bar{x}^i:=(\bar{x}^i_0,\dots, \bar{x}^i_{n})$ where 
$\bar{x}_j^i:={\phi^i}^*x_j^i$ and $x^i_j$ is the function defined in Subsection \ref{einsforg} for $g_i$ with boundary values $x_j^1|_{\rho_1=0}=x_j^2|_{\rho_2=0}\in C_0^\infty(\Gamma)$ for all $j$, 
Now set $\psi:U\to \psi(U)\subset [0,\eps)\x\Gamma$ such that 
$\bar{x}_j^1=\psi^*\bar{x}_j^2.$  This is a diffeomorphism near $p$ and  
moreover ${\phi^1}^*g_1$ and ${\phi^2}^*g_2$ coincide to infinite order at $t=0$ by Lemma \ref{infiniteorder}, 
so the coordinates $\bar{x}^1$ and $\bar{x}^2$ satisfy $\Delta_{{\phi^1}^*g_1}(\bar{x}_j^1-\bar{x}_j^2)=O(t^{\infty})$ 
for all $j$. Since  $\bar{x}^1,\bar{x}^2$ have the same boundary values, they agree to order $O(t^n)$ using the construction
of $F_{n-1}$ in \eqref{constfj}. But since $\mc{S}_1(x_j^1|_\Gamma)=\mc{S}_2(x^2_j|_{\Gamma})$ on $\Gamma$, one has 
$\bar{x}_j^1=\bar{x}^2_j+O(t^{n+1})$ near $p$, which again by induction and \eqref{indiciel} shows that $\bar{x}_j^1=\bar{x}^2_j+O(t^{\infty})$ near $p$.
In particular, setting  $\bar{g}_1:=(\bar{x}^1_0)^2{\phi^1}^*g_1$, $\bar{g}_2:=\psi^*((\bar{x}^2_0)^2{\phi^2}^*g_2)=(\bar{x}_0^1)^2\psi^*({\phi^2}^*g_2)$, one obtains that $\bar{g}_1=\bar{g}_2+O((\bar{x}^1_0)^{\infty})$, i.e. the metrics 
agree to infinite order in the coordinates $\bar{x}^1$.
Thus from Lemma \ref{system2}, $\bar{g}_1$ and $\bar{g}_2$ both
satisfy the same system \eqref{einsteing} and agree to infinite order at the boundary $\{t=0\}$ near $p$ in the coordinate system $\bar{x}^1$, 
so by Proposition \ref{unique1}, we deduce that ${\phi^1}^*g_1=\psi^*{\phi^2}^*g_2$ and this ends the proof by setting $F:=\phi^2\circ\psi\circ(\phi^1)^{-1}$. 
\qed\\

We finish by the following Corollary, similar to Corollary \ref{corol1}. 
\begin{cor}\label{corol2}
Let $G_i(z,z')$ be the Green kernel for $g_i$, $i=1,2$. Then $\mc{S}_1|_{\Gamma}=\mc{S}_2|_{\Gamma}$ implies that there exists
an open set $U_1$ such that $G_1(z,z')=G_2(F(z),F(z'))$ for all $(z,z')\in (U_1\x U_1)\setminus\{z=z'\}$.
\end{cor}
\noindent \textsl{Proof}:  We first take $y'\in U_1\cap \Gamma$, and consider $T_1(z):=[\rho_1(z')^{-n}G_1(z,z')]_{z'=y'}$ 
and $T_2(z):=[F^*\rho_2(z')^{-n}G_2(F(z),F(z'))]|_{z'=y'}$. They both satisfy $\Delta_{g_1}T_i(z)=0$
for $z\in U_1$ and by Lemma \ref{noyausc} and the assumption $\mc{S}_1|_{\Gamma}=\mc{S}_2|_{\Gamma}$,
we have that $T_1-T_2=O(\rho_1^{n+1})$ near $\Gamma\setminus\{y'\}$, so by induction on \eqref{indiciel},
$T_1=T_2+O(\rho^1_\infty)$ in $U_1\setminus\{y'\}$, and then by the unique continuation result of Mazzeo \cite{Ma},
$T_1=T_2$ in the same set. Now this means that for $z'\in U_1$, $z\to G_1(z',z)$ and $z\to G_2(F(z'),F(z))$ are
harmonic for $\Delta_{g_1}$ in $U_1\setminus\{z'\}$, and they coincide to order $\rho_1^{n+1}$ at $\Gamma$, 
so again by unique continuation they are equal.
\qed

\subsection{Proof of Theorem \ref{th3}}

Using Corollary \ref{corol2} and the fact that $(X_1,g_1)$ and $(X_2,g_2)$ Einstein, and by Theorem 5.2 of \cite{DTK} are analytic in harmonic coordinates, 
it suffices to apply Theorem 4.1 of \cite{LUT}, which is essentially the same as Proposition \ref{lut1} but for a complete manifold.

\end{document}